\documentclass[12pt,a4paper,reqno]{amsart}
\usepackage{amsmath,amsthm}
\usepackage{amsfonts}
\usepackage{amssymb}
\allowdisplaybreaks

\newcommand{\be}{\begin{equation}}
\newcommand{\ef}{\end{equation}}
\chardef\bslash=`\\ % p. 424, TeXbook
\newtheorem*{thm*}{Theorem}

\theoremstyle{definition}

\newtheorem*{remark*}{Remarks}
\newtheorem*{defn*}{Definition}

\theoremstyle{remark}
\newcommand{\G}{\Gamma}
\newcommand{\wt}{\widetilde}
\newcommand{\wh}{\widehat}
\newcommand{\fc}{\frac}

\newcommand{\iy}{\infty}
 \renewcommand{\sectionmark}[1]{}

\renewcommand{\Re}{\operatorname{Re}}

\newcommand{\ve}{\varepsilon}

\newcommand{\Te} {Teichm\"{u}ller}

\newcommand{\const}{\operatorname{const}}

 \usepackage{amsfonts}
\newcommand{\field}[1]{\mathbb{#1}}
\newcommand{\g}{\gamma}

\newcommand{\dl}{\delta}
\newcommand{\D}{\field{D}}
\newcommand{\om}{\omega}
\newcommand{\z}{\zeta}
\newcommand{\ov}{\overline}
\newcommand{\vp}{\varphi}
\newcommand{\hC}{\wh{\field{C}}}
\newcommand{\C}{\field{C}}
\newcommand{\R}{\field{R}}

\newcommand{\B}{\mathbf{B}}
\newcommand{\T}{\mathbf{T}}
\newcommand{\Belt}{\operatorname{Belt}}

\newcommand{\Fib}{\operatorname{Fib}}
\newcommand{\Teich}{\operatorname{Teich}}

\newcommand{\vk} {\varkappa}
\newcommand{\x} {\mathbf x}
\renewcommand{\a} {\alpha}

\newcommand{\ld}{\lambda}
\newcommand{\kp}{\kappa}

\begin{document}

\title{Quasiconformal deformations preserving Hilbert's norm and their applications}
\author{Samuel L. Krushkal}

\begin{abstract} This paper focuses on estimating the Taylor coefficients for Hilbert spaces of holomorphic functions  on the disk using intrinsic features of univalent functions and of Teichm\"{u}ller spaces.
Estimating these coefficients has a long history but still remains an important problem in many geometric and physical applications of complex analysis.

We construct quasiconformal deformations of holomorphic functions preserving their Hilbert norm.
Such deformations play a crucial role in this subject.
Among their applications, a rather complete solution of an old Hummel-Scheinberg-Zalcman problem
is obtained.
\end{abstract}

\date{\today\hskip4mm({estHilb.tex})}

\maketitle

\bigskip

{\small {\textbf {2020 Mathematics Subject Classification:} Primary: 30C50, 30C62, 30E10, 30F60;
Secondary 30C75, 31A05, 32L05, 32Q45}

\medskip

}

\bigskip

\markboth{Samuel L. Krushkal}{Quasiconformal deformations preserving Hilbert's norm and their applications} \pagestyle{headings}

\bigskip\bigskip
\centerline{\bf 1. INTRODUCTION}

\bigskip
One of new interesting phenomena in geometric complex analysis discovered recently is a natural connection of distortion results for the general classes of holomorphic functions $f$ from some Banach spaces with the appropriate collections of univalent functions and Teichm\"{u}ller spaces.

A crucial role here is played by special quasiconformal deformations $h$  of such  functions; these deformations must preserve a finite set of distinguished Taylor coefficients of $f$ and the norm of $f$.

The existence of such deformations was established in \cite{Kr2} and extended in \cite{Kr6} to spaces
of holomorphic functions with weighted $L_p$-norm. Such deformations found various deep applications
though the indicated classes are rather narrow. Their functions  are holomorphic on $\D$ and have the
growth
$$
f(z) = O(1/(1 - |z|)^p) \quad \text{as} \ \ |z| \to 1
$$
with $0 \le p < 2$. The corresponding solutions of the Schwarz differential equation
$$
(w^{\prime\prime}(z)/w^\prime(z))^\prime - (w^{\prime\prime}(z)/w^\prime(z))^2/2 = f(z)
$$
are asymptotically conformal on the boundary; hence, such $f$ can run only over a nonwhere dense subsets in the space of all Schwarzians.

So there is an interesting open problem to find new spaces of holomorphic functions on the disk
(more generally, on arbitrary bounded quasidisks) admitting quasiconformal perturbations which vary
the prescribed coefficients and preserve the norm.

In this paper, we solve this problem for a wide collection of Hilbert spaces and apply this to solving
for such spaces the well-known old Hummel-Scheinberg-Zalcman problem  \cite{HSZ} on estimating the coefficients of nonvanishing (zero free) holomorphic functions in Banach spaces.

\bigskip\bigskip
\centerline{\bf 2. QUASICONFORMAL DEFORMATIONS PRESERVING THE HILBERT}
\centerline{\bf NORM}

\bigskip
Let $H$ be a complex Hilbert space of holomorphic functions
$$
f(z) = c_0 + c_1 z + \dots + c_n z^n + \dots
$$
in the unit disk $\D = \{|z| < 1\}$ with scalar product $(f, g)_X$ containing the linear subspace spanned by the system $\{z^n\}_0^\iy$.

We also assume that if $f \in H$ and the complementary set $\C \setminus \ov{f(\D)}$   is a domain, then
the compositions $h \circ f$  with $h$ holomorphic on $\ov{f(\D)}$ also belong to $H$. In fact, we need this only for fractions
$$
c/(\z - f(z)) \quad \text{with} \ \ \z \in \C \setminus \ov{f(\D)}, \ \ c = \const. 
$$

For such spaces, we have the following theorem which underlies other results of this paper.

\bigskip
Fix a natural $n >  1$ and consider the collections $\mathbf d = (d_0, d_1, \dots, d_n) \in \C^{n+1}$ with the Euclidean norm $|\mathbf d|$.

\bigskip\noindent
{\bf Theorem 1}. {\it If a function $f(z) = c_j z^j + c_{j+1} z^{j+1} + \dots \in H$ with $c_j \ne 0, \ 0 \le j < n$, is not a polynomial of degree $n_1 \le n$ and $f(\D)$ does not contain the points from a distinguished (sufficiently distinct from the origin) disk $\D(w_0, r) = \{|w - w_0| < r\}$ , then there exists a positive number $\ve_0$ such that for every point
$$
{\mathbf d} = (d_{j+1}, \dots, \ , d_n) \in \C^{n-j}
$$
and every $a \in \R$ with
$$
|\mathbf  d| \le \ve, \ \ |a| \le \ve, \ \ \ve < \ve_0,
$$
there is a quasiconformal automorphism $h$ of the extended complex plane $\hC$, which is conformal
outside of the disk $\D(w_0, r)$ and satisfies the conditions:

\noindent
(a) $h^{(k)}(f(0)) = k! d_k = k! (c_k + d_k), \  k = j + 1, \dots, n$;

\noindent
(b) $\|h \circ f\|_H =  \|f \|_H + a$.

\noindent
(c) The Beltrami coefficient $\mu_h$ of $h$ satisfies $\|\mu_h \|_\iy \le M \ve$.

\noindent
The quantities $\ve_0$ and $M$ depend only on $f$ and $n$}.

\bigskip
Recall that the Beltrami equation is of the form
\be\label{1}
\partial_{\ov w} h = \mu(w) \partial_w h,
\end{equation}
where the partial derivatives $\partial_{\ov w}, \ \partial_w h$ are understand generically as distributional. Quasiconformal homeomorphisms are the generalized solutions of this equation.

Clearly, this theorem holds for bounded (with finite $L_\iy$-norm) functions from $H$.
Quasiconformal deformations of such functions given by Theorem 1 preserve their $H$-norm, but generically increase $L_\iy$-norm.

The special cases of this theorem have been proved and applied in \cite{Kr2}, \cite{Kr6}. In these cases,
the system of linear algebraic equations arising in the proof has diagonal matrix, which is caused by the mutual orthogonality of powers $z^n, \ n = 0, 1, \dots$, on the circles $\{|z| = R\}$.

\bigskip\bigskip
\centerline{\bf 3. PROOF OF THEOREM 1}

\bigskip
We apply the standard representation of quasiconformal homeomorphisms of $\hC$ using the integral operators
$$
T\rho = - \frac{1}{\pi} \iint\limits_{\D(w_0, r)} \frac{\rho(\zeta) d \xi d \eta}{\z - w}, \quad
\Pi \rho = \partial_w T\rho = - \frac{1}{\pi} \iint\limits_{\D(w_0, r)}  \frac{\rho(\zeta) d \xi d \eta}{(\z - w)^2}
$$
on the space $L_p(\D(w_0, r)), \ p > 2$ (the second integral exists as a principal Cauchy value). The required quasiconformal automorphism $h = h^\mu$ has the form
 \be\label{2}
h(w) = w - \frac{1}{\pi} \iint\limits_{\D(w_0, r)}  \frac{\rho(\zeta) d \xi d \eta}{\zeta - w} = w + T \rho(w);
\end{equation}
its Beltrami coefficient $\mu = \mu_h$ (extended by zero outside of $\D(w_0, r)$) generates after substitution of (2) into the equation (1) the series
$$
\rho = \mu + \mu \Pi \mu + \mu \Pi (\mu\Pi \mu) + ... \ .
$$
If $\|\mu \|_\iy < \kp < 1$, then the last series is convergent in $L_p(\C)$ for some $p = p(\kp) > 2$.
The well-known properties of operators $T$ and $\Pi$ imply for any disk $\D_R = \{|w| < R\} \ (0 < R < \iy)$ the estimates
$$
\|\rho\|_{L_p(\D_R)},  \ \ \| \Pi \rho \|_{L_p(\D_R)} \le  M_1(\kp, R, p) \|\mu\|_{L_\iy(\C)};
\|h\|_{C(\D_R)} \le  M_1(\kp, R, p) \|\mu \|_\iy
$$
(valid for a fixed disk $\D(w_0, r)$). Hence,
 \be\label{3}
h(w) = w + T \mu(w) + \om(w)
$$
with $\|\om \|_{C(\D_R)} \le M_2(\kp, R, p) \|\mu \|_\iy^2$.
The bounds $M_1, \ M_2 \to 0$ as $\kp \to 0$ or $r \to 0$.

Using the pairing
$$
\langle\nu,\vp\rangle = - \fc{1}{\pi} \iint\limits_{\D(w_0, r)} \nu(\zeta) \vp(\zeta) d \xi d \eta,
\quad \nu \in L_\iy(\D(w_0, r)), \ \ \vp \in L_1(\D(w_0, r)),
\end{equation}
one can rewrite the above representation in the form
 \be\label{4}
h(w) = w + \sum\limits_1^\iy \langle\mu, \vp_k\rangle \ (w - c_0^0)^k + \om(w), \quad
\vp_k(\z) = \fc{1}{(\z - c_0^0)^k},
\end{equation}
where
$$
\vp_k(\z) = 1/(\z - c_0)^k, \quad k = j + 1, ... \ , n.
$$
The representation (4) with $k = j + 1, \dots, n$ and condition $(a)$ of the theorem provide the first group of equalities to determine the desired Beltrami coefficient $\mu$:
 \be\label{5}
k! d_k = \langle \mu,\vp_k \rangle + \om^{(k)}(c_0) = \langle\mu, \vp_k\rangle + O(\|\mu \|_\iy^2), \quad k = j + 1, ... \ , n.
\end{equation}

Now we apply the remained initial function $\vp_j$ to vary $\|f\|_X$. This variation is constructed as follows. First, noting that
$$
- \fc{1}{\pi} \iint\limits_{|\z-w_0|<1} \fc{d \xi d \eta}{\z - w} =
\begin{cases} \ov{w} - \ov{w_0} &  \text{if} \quad |w - w_0| < 1,  \cr
               1/(w - w_0)   &  \text{if} \quad |w - w_0| > 1,
\end{cases}
$$
one obtains for functions $\mu \in C^2(\ov{\D(w_0, r)})$ and $|w - w_0| \ge r_0 > r$ the asymptotic equality
 \be\label{6}
- \fc{1}{\pi} \iint\limits_{\D(w_0, r)} \fc{\mu(\z) d \xi d \eta}{\z - w} =
\fc{\mu(w_0)}{w - w_0} + O(r)
\end{equation}
with uniform estimate of the reminder. Hence, for $w \in \D(w_0, r)$ with a small $r > 0$,
$$
h(w) = \mu(w_0) (\ov w - \ov w_0) + O_r(\|\mu \|_\iy^2)
$$
(this is an affine map up to a quantity of higher order), while the equalities (5) can be rewritten in the
form
 \be\label{7}
k! d_k = \fc{\mu(w_0)}{(c_0 - w_0)^{k+1}} + O_r(\|\mu \|_\iy^2) + O(r) \quad k = j + 1, ... \ , n.
\end{equation}
Note that the reminder terms in (6), (7) tend to $0$ as $r \to 0$. Together with (3), this yields
$$
\|h \circ f\|_H^2 = (h \circ f, \ov{h \circ f})_H
= (f + T\mu \circ f, \ov{f + T\mu \circ f})_H
= \|f\|_H^2 + 2 \Re (\ov f, T\mu \circ f)_H + O_r(\|\mu\|_\iy^2).
$$
In view of (6), this can be written for sufficiently large $|w_0|$ and small $r$ in the form
 \be\label{8}
\|h \circ f\|_H^2 - \|f\|_H^2 = 2 \Re (f, \ov{\mu(w_0)/(w - w_0)} \circ f)_H + O(\|\mu\|_\iy^2).
\end{equation}
The first term in the right-hand side of (8) is well defined on the unit disk.

Noting that, in view of the assumption $f(0) \ne w_0$, the function
$$
\phi(w) = 1/(w - w_0)
$$
is distinct from a linear combination of the fractions $\vp_{j+1}, ..., \vp_n$, one can define on the
linear span of $\vp_{j+1}, ..., \vp_n, \phi$ the linear functional $l(\vp)$ orthogonal to all
$\vp_{j+1}, ..., \vp_n$ and with $l(\phi) = 1$, which defines (after Hahn-Banach extension to
$A_1(\D(w_0, r))$ and then to $L_1(\D(w_0, r))$ ) the corresponding Beltrami coefficient $\mu_0$ supported on $\D(w_0, r)$ with
$$
\langle \mu_0, \vp_k\rangle = 0, \ \ k = j + 1, \dots, n; \ \ \langle \mu_0, \phi\rangle =
\|\mu_0\|_\iy = 1.
$$
Now we seek the desired Beltrami coefficient $\mu$ in the form
  \be\label{9}
\mu(w) = \sum\limits_{j+1}^n \xi_k \ov \vp_k + \tau \mu_0
\end{equation}
(extended by zero outside of the disk $\D(w_0,r)$) with unknown constants
$\xi_{j+1}, ... \ , \xi_n, \tau$
to be determined from the equalities (7) and (8).

Substitution of (9) into (7) and (8) provides the nonlinear equations
 \be\label{10}
\begin{aligned}
k! d_k &= \sum\limits_{l=j+1}^n \xi_l \langle \ov \vp_l, \vk\rangle_H + O(\|\mu \|_\iy^2), \quad
k = j + 1, ... \ , n, \\
\|h \circ f_0 \|_H^2 - \|f_0 \|_H^2 &= \tau \mu_0 + O(\|\mu \|_\iy^2)
\end{aligned}
\end{equation}
for determining $\xi_k$ and $\tau$. This equations define a nonlinear holomorphic map
$$
\mathbf y = W(\mathbf x) = W^\prime(\mathbf 0) \mathbf x + O(|\mathbf x|^2),
$$
of a small neighborhood $U_0$ of the origin in $\C^{(n - j) + 1}$. Its linearization
$\mathbf y = W^\prime(\mathbf 0) \mathbf x$ provides a linear map of $\C^{n -j + 1}$ whose
Jacobian is a nonzero Gram determinant generated by the linearly independent functions
$\vp_{j+1}, \dots, \vp_n, \phi$.
Therefore, $\mathbf x \mapsto W^\prime(\mathbf 0) \mathbf x$ is a linear isomorphism of the space
$\C^{n - j + 1}$ onto itself, and one can apply to $W$ the inverse mapping theorem, which provides
the assertion of Theorem 1.

\bigskip\bigskip
\centerline{\bf 4. APPLICATIONS}

\bigskip\noindent
{\bf 4.1. Hyperbolically bounded holomorphic functions}.
The first application of Theorem 1 concerns a new method for estimating the coefficients of holomorphic functions introduced in \cite{Kr6} and extends this method to functions from a rich collection of Hilbert spaces. Though estimating of coefficients of is a very old problem of complex analysis, it still remains important having a lot of geometric and physical applications.
The indicated approach is based on deep features of univalent functions and on Teichm\"{u}ller space theory.

First consider the Hilbert spaces $H$ whose functions
$f(z) = c_0 + c_1 z + \dots + c_n z^n + \dots, \ z \in \D$ have the growth
$$
f(z) = O(1/(1 - |z|)^2) \quad \text{as} \ \ |z| \to 1.
$$
For our goals it is more conveniently deal with the  quadratic differentials $f(z) dz^2$.
Such functions form the complex Banach space
$\B_2$ with norm
$$
\|f\|_{\B_2} = \sup_\D (1 - |z|^2)^2 |f(z)|,
$$
which is dual to the subspace $A_1(\D)$ of integrable holomorphic functions (moreover, the mean value inequality yields for $f \in A_1(\D)$ that $\|f\|_{\B_2} \le \|f\|_{A_1}$; hence $A_1(\D) \subset\B_2$).

In this section we consider the Hilbert spaces $H$ embedded in $\B_2$ and also assume that
 \be\label{11}
\|f\|_{\B_2} \le \const \|f\|_H  
\end{equation}
with a constant not depending on $f \in H$. 

The space $\B_2$ has a special interest because its points are the {\bf Schwarzian derivatives}
$$
S_w(z) = \left(\frac{w^{\prime\prime}(z)}{w^\prime(z)}\right)^\prime
- \frac{1}{2} \left(\frac{w^{\prime\prime}(z)}{w^\prime(z)}\right)^2
$$
of locally univalent holomorphic functions $w(z) = a_0 + a_1 z+ \dots$ on $\D$ which are determined by
$f \in \B_2$ as solutions of the differential equation
 \be\label{12}
(w^{\prime\prime}(z)/w^\prime(z))^\prime - (w^{\prime\prime}(z)/w^\prime(z))^2/2 = f(z)
\end{equation}
satisfying the prescribed initial conditions
 \be\label{13}
w(0) = a_0, \ \ w^\prime(0) = a_1, \ \ w^{\prime\prime}(0) = a_2.
\end{equation}
Equivalently, $w$ is the ratio $\eta_1/\eta_2$ of two linearly independent solutions $\eta_1, \eta_2$ of the linear differential equation
  \be\label{14}
2 \eta^{\prime\prime}(z) + f(z) \ \eta(z) = 0
\end{equation}
with $\eta_1(0) = 0, \ \eta_1^\prime(0) = 1$ and $\eta_1(0) = 1, \ \eta_1^\prime(0) = 0$.

The coefficients $c_n$ are uniquely determined by coefficients $a_m$ of the corresponding $w(z)$,
and the problem is to find the explicit bounds for $c_n$ having the bounds for $a_m$.

The Schwarzians $S_w$ of univalent solutions $w(z)$ with quasiconformal extension to the complementary disk $\D^* = \{z \in \hC: \ |z| > 1\}$ fill in $\B_2$ a bounded domain modeling
the universal Teichm\"{u}ller space $\T$.

\bigskip
The main result of this section is the following theorem solving implicitly the indicated
problem. It essentially extends the results established in \cite{Kr6} for spaces with $L_p$-norm.

\bigskip\noindent
{\bf Theorem 2}. {\it (a) Let $G$ be a complex submanifold $G$ of a Hilbert space $H$, which satisfies (11), and let the boundary of $G$ have the common points with $\partial \T$. Then the function $f_0(z) = c_0 + c_1 z + \dots$ maximizing the first coefficient $c_1$ on $\ov G = G \cup \partial G$ is the Schwarzian $S_{w_0}$ of a univalent function $w_0(z)$ with maximal second coefficient $a_2$ among all univalent functions $w$ on the disk $\D$ whose Schwarzians $S_w$ are the points of $\ov G$. }

\noindent
{\it (b) The coefficients $a_m$ and non-zero $c_n$ are estimated for all $m \ge 3$ and $n \ge2$ by
 \be\label{15}
|a_m| \le |a_m^0|
\end{equation}
and
 \be\label{16}
|c_n| \le \max(|c_1^0|, |c_n^0|).
\end{equation}
}

In particular, if $|c_n^0| \le |c_1^0|$ for all $n \ge 2$, then all $c_n$ are bounded by this $|c_1^0|$.

\bigskip
Using Theorem 1 one can prove Theorem 2 in the lines of \cite{Kr6}, embedding by (11) the given  space $H$ into the universal Teichm\"{u}ller space $\T$ and using Bers' fiber space $\Fib (\T)$ over $\T$.
In view of importance of Theorem 2, we briefly outline this proof in the next section.

This theorem allows one to reduce estimating the coefficients of arbitrary holomorphic functions
to corresponding classes of univalent functions whose theory is widely developed and provides many
powerful methods. A crucial point here is that Hilbert's norm does not be increased.

\bigskip\noindent
{\bf 4.2. Functions with higher boundary singularities}. Now fix an integer $p > 2$ and consider the Hilbert spaces $H$ whose functions have the growth
$$
f(z) = O(1/(1 - |z|)^p) \quad \text{as} \ \ |z| \to 1.
$$
All such holomorphic functions on $\D$ form the Banach space $\B_p(D)$ with norm
$$
\|\vp\|_{\B_p} = \sup_\D (1 - |z|^2)^p |\vp(z)|;
$$
it is dual to the space $A_p(\D)$ of integrable holomorphic functions on the disk $\D$ with the weight
$(1 - |z|^2)^{2-p}$.

The relation between the spaces $\B_p$ and $\B_2$ is characterized by the following theorem which is a special case of the approximation results established in \cite{Kr4}.

\bigskip\noindent
{\bf Theorem 3}. {\it For any function $f \in \B_p(\D)$ there exists a sequence of rational functions
with poles of order two on the unit circle $\mathbb S^1 = \{|z| = 1\}$ of the form
 \be\label{17}
r_n(z) = \sum\limits_1^n \fc{d_j}{(z - a_j)^2}, \quad
\sum\limits_1^n |d_j| > 0,
\end{equation}
such that $\lim\limits_{n\to \iy} \|r_n - \vp\|_{\B_{p+1}(\D)} = 0$.}

In fact, such approximation is valid for all bounded simply connected domains with quasiconformal boundaries (quasidisks) and in the case of the disk $\D$ (or half-plane) the coefficients $d_j$ can be chosen to be real.

Note also that the set of hyperbolically boundended functions $\vp$ in the corresponding space $\B_2(D)$ for a such domain $D$ approximated in $\B_2$-norm by rational functions with poles of order two on the boundary $\partial D$ is nonwhere dense in this space.
\footnote{Recall that $\|\vp\|_{\B_2(D)} = \sup_D \ld_D(z)^{-2} |\vp(z)|$, where $\ld_D(z) |dz|$ is
the hyperbolic metric of $D$ of Gaussian curvature $- 4$.}

Theorems 2 and 3 provide together the approximative bounds for coefficients $c_n$ of functions $f \in \B_p$.

\bigskip\noindent
{\bf 4.3. Nonvanishing functions. The Hummel-Scheinberg-Zalcman problem}.
Such functions arise and play an essential role in many fields of complex analysis, in particular,
as the derivatives of locally univalent functions and of covering maps and have been investigated by
many authors. There were several deep conjectures about such functions including the eminent Krzyz and Hummel-Scheinberg-Zalcman conjectures for nonvanishing holomorphic functions from the Hardy spaces
$H^p, \ 1< p \le \iy$ (see, e.g., \cite{HSZ}, \cite{Kr5}, \cite{Kz} and the references cited there).

The paper \cite{HSZ} also contains a collection of interesting open questions on estimating the coefficients in various classes of nonvanishng holomorphic functions on the disk, among those the following general problem:

\bigskip\noindent
{\it Find the value of
 \be\label{19}
C_n(B) = \sup_{\mathcal B(B)} |c_n|,
\end{equation}
where $\mathcal B(B)$ denotes the class of functions $f(z) = \sum_0^\iy c_n z^n, \ |z| < 1$, which belong
to the (complex) Banach space $B$ and satisfy} $\|f\|_B \le 1, \ f(z) \ne 0$.

\bigskip
All such spaces are compact in the weak topology of locally uniform convergence on the disk $\D$, which insures the existence of extremal functions.

\bigskip
The existence of variations given by Theorem 1 provides a possibility to solve the
Hummel-Scheinberg-Zalcman problem for Hilbert spaces $H$ satisfying some natural additional conditions: $H$ must admit the assumption (11) (be embedded into $\B_2$) and contain together with any function
$f(z)$ its homotopy functions $f_r(z) = f(r z), \ 0 \le r \le 1$, defining a curve in $H$  connecting
$f$ with the origin $f(z) \equiv 0$ of $H$.

\bigskip\noindent
{\bf Theorem 4}. {\it For any Hilbert space $H$ satisfying the indicated assumptions, the function $f_0(z)$ maximizing the first coefficient $c_1$ on $\mathcal B(H)$ is the Schwarzian of a univalent solution $w_0$ of the corresponding equation (12) with maximally admissible value $|a_2|$ in (13),
and the larger coefficients $c_n, \ n > 1$, of every function $f \in \mathcal B(H)$ are estimated
via (16) by coefficients $c_1^0, \ c_n^0$ of $S_{w_0}$.
}

\bigskip\bigskip
\centerline{\bf 5. REMARKS ON TEICHM\"{U}LLER SPACES}

\bigskip
As was mentioned above, the proof of Theorems 2 and 4 essentially relies on the deep geometric and
analytic features of Teichm\"{u}ller spaces. Thus we briefly recall the needed results from this theory;
the details can be found, for example, in \cite{Be}, \cite{GL}, \cite{Le}.
It is technically more convenient to deal with functions from $\Sigma_Q$.

\bigskip\noindent
{\bf 5.1}.
The {\bf universal Teichm\"{u}ller space} $\T = \Teich(\D)$ is the space of quasisymmetric homeomorphisms of the unit circle $\mathbb S^1$ factorized by M\"{o}bius maps;  all Teichm\"{u}ller spaces have their biholomorphic copies in $\T$.

The canonical complex Banach structure on $\T$ is defined by factorization of the ball of the Beltrami coefficients (or complex dilatations)
$$
\Belt(\D)_1 = \{\mu \in L_\iy(\C): \ \mu|\D^* = 0, \ \|\mu\| < 1\},
$$
letting $\mu_1, \mu_2 \in \Belt(\D)_1$ be equivalent if the corresponding  quasiconformal maps $w^{\mu_1}, w^{\mu_2}$ (solutions to the Beltrami equation $\partial_{\ov{z}} w = \mu \partial_z w$
with $\mu = \mu_1, \mu_2$) coincide on the unit circle $\mathbb S^1 = \partial \D^*$ (hence, on $\ov{\D^*}$). Such $\mu$ and the corresponding maps $w^\mu$ are called $\T$-{\it equivalent}.

The following important lemma from \cite{Kr5} allows one to use some other normalizations of
quasiconformally extendable functions.

\bigskip\noindent
{\bf Lemma 1}. {\it For any Beltrami coefficient $\mu \in \Belt(\D^*)_1$ and any $\theta_0 \in [0, 2 \pi]$, there exists a point $z_0 = e^{i \a}$ located on $\mathbb S^1$ so that
$|e^{i \theta_0} - e^{i \a}| < 1$ and such that for any $\theta$ satisfying
$|e^{i \theta} - e^{i \a}| < 1$ the equation
$\partial_{\ov z} w =  \mu(z) \partial_z w$
has a unique homeomorphic solution $w = w^\mu(z)$, which is holomorphic on the unit disk $\D$
and satisfies
 \be\label{19)}
w(0) = 0, \quad w^\prime(0) = e^{i \theta}, \quad w(z_0) = z_0.
\end{equation}
Hence, $w^\mu(z)$ is conformal and does not have a pole in $\D$ \ (so
$w^\mu(z_{*}) = \iy$ at some point $z_{*}$ with $|z_{*}| \ge 1$).  }

\bigskip
In particular, this lemma allows one to define the Teichm\"{u}ller spaces using the quasiconformally extendible  univalent functions $w(z)$ in the unit disk $\D$ normalizing these functions by
$$
w(0) = 0, \quad w^\prime(0) = e^{i \theta}, \quad w(1) = 1.
$$
All such functions are holomorphic in the disk $\D$.

\bigskip\noindent
{\bf 5.2}.
The Teichm\"{u}ller space $\T_1 = \Teich(\D_{*})$ {\bf of the punctured disk} $\D_{*} = \D \setminus \{0\}$ is formed by classes $[\mu]_{\T_1}$ of $\T_1$-{\bf equivalent} Beltrami coefficients $\mu \in \Belt(\D)_1$ so that the corresponding quasiconformal automorphisms $w^\mu$ of the unit disk coincide on both boundary components (unit circle $\mathbb S^1$ and the puncture $z = 0$) and are homotopic on $\D \setminus \{0\}$.
This space can be endowed with a canonical complex structure of a complex Banach manifold
and embedded into $\T$ using uniformization of $\D_{*}$ by a cyclic parabolic Fuchsian
group acting discontinuously on $\D$ and $\D^*$. The functions $\mu \in L_\iy(\D)$ are lifted to
$\D$ as the Beltrami measurable $(-1, 1)$-forms  $\wt \mu d\ov{z}/dz$ in $\D$ with respect to
$\G$, i.e., via $(\wt \mu \circ \g) \ov{\g^\prime}/\g^\prime = \wt \mu, \ \g \in \G$,
forming the Banach space $L_\iy(\D, \G)$; we extend these $\wt \mu$ by zero to $\D^*$. Then $\T_1$ is canonically isomorphic to the subspace $\T(\G) = \T \cap \B_2(\G)$,
where $\B_2(\G)$ consists of elements $\vp \in \B_2$ satisfying $(\vp \circ \g) (\g^\prime)^2 = \vp$
in $\D^*$ for all $\g \in \G$.

Due to {\bf the Bers isomorphism theorem}, {\it the space $\T_1$ is biholomorphically isomorphic
to the Bers fiber space
$$
\Fib(\T) = \{(\phi_\T(\mu), z) \in \T \times \C: \ \mu \in \Belt(\D)_1, \ z \in w^\mu(\D)\}
$$
over the universal Teichm\"{u}ller space $\T$ with holomorphic projection $\pi(\psi, z) = \psi$} (see \cite{Be}).

This fiber space is a bounded hyperbolic domain in $\B_2 \times \C$ and represents the collection of domains $D_\mu = w^\mu(\D)$ as a holomorphic family over the space $\T$. For every $z \in \D$,  its
orbit $w^\mu(z)$ in $\T_1$ is a holomorphic curve over $\T$.

The indicated isomorphism between $\T_1$ and $\Fib(\T)$ is induced by the inclusion map \linebreak
$j: \ \D_{*} \hookrightarrow \D$ forgetting the puncture at the origin via
 \be\label{20}
\mu \mapsto (S_{w^{\mu_1}}, w^{\mu_1}(0)) \quad \text{with} \ \
\mu_1 = j_{*} \mu := (\mu \circ j_0) \ov{j_0^\prime}/j_0^\prime,
\end{equation}
where $j_0$ is the lift of $j$ to $\D$.

The Bers theorem is valid for Teichm\"{u}ller spaces $\T(X_0 \setminus \{x_0\})$ of all punctured hyperbolic Riemann surfaces $X_0 \setminus \{x_0\}$; we use only its special case.

\bigskip\noindent
{\bf 5.3}.
The spaces $\T$ and $\T_1$ can be weakly (in the topology generated by the spherical metric on $\hC$) approximate by finite dimensional Teichm\"{u}ller spaces $\T(0, n)$ of punctured spheres (Riemann surfaces of genus zero)
$$
X_{\mathbf z} = \hC \setminus \{0, 1, z_1 \dots, z_{n-3}, \iy\}
$$
defined by ordered $n$-tuples $\mathbf z = (0, 1, z_1, \dots, z_{n-3}, \iy), \ n > 4$ with distinct
$z_j \in \C \setminus \{0, 1\}$ (the details see, e.g., in \cite{Kr3}).

Another canonical model of $\T(0, n) $ is obtained again using the uniformization. This space
is biholomorphic to a bounded domain in the complex Euclidean space $\C^{n-3}$.

Note also that all Teichm\"{u}ller spaces are complete metric spaces with intrinsic Teichm\"{u}ller metric defined by quasiconformal maps. By the Royden-Gardiner theorem, this metric equals the hyperbolic  Kobayashi metric determined by the complex structure (see, e.g., \cite{EKK}, \cite{GL}, \cite{Ro}).

\bigskip\bigskip
\centerline{\bf 6. PROOF OF THEOREM 2}

\bigskip
The proof given below follows the lines of \cite{Kr3}, \cite{Kr5}. In view importance of this theorem, we outline the main steps of this proof. It involves lifting the functional $J(w) = c_n$ onto the spaces $\T$ and $\T_1$.

\noindent
{\it(a)} \ Consider the classes $S_{z_0, \theta}(D)$ of univalent functions in $\D$ with expansions
$$
f(z) = a_1 z^{-1} + a_2 z^2 + \dots
$$
with
$a_1 = e^{i \theta}, \ - \pi < \theta \le \pi$, having the fix point at $z_0 \in \mathbb S^1$ and
admitting quasiconformal extensions to $\hC$ and take their union
$$
\wh S = \bigcup_{z_0 \in \mathbb S^1, \theta \in (-\pi,\pi]} S_{z_0,\theta}.
$$
The Beltrami coefficients $\mu_f(z) = \partial_{\ov z} f/\partial_z f$ of these extensions run over the unit ball
$$
\Belt(D^*)_1 = \{\mu \in L_\iy(\C): \ \mu(z)|D = 0, \ \ \|\mu\|_\iy  < 1\}.
$$

Now pass to the inverted functions $F_f(z) = 1/f(1/z)$ which form the corresponding classes $\Sigma_{z_0,\theta}$ of nonvanishing univalent functions on the disk $\D^*$ with expansions
$$
F(z) =  e^{- i \theta} z + b_0 + b_1 z^{-1} + b_2 z^{-2} + \dots, \quad  F(z_0) = z_0,
$$
and let
$\Sigma^0 = \bigcup_{z_0,\theta} \Sigma_{z_0,\theta}$.

The coefficients $a_n$ of $f(z)$ and the corresponding coefficients $b_j$ of $F_f(z)$ are related by
$$
b_0 + e^{2i \theta} a_2 = 0, \quad b_n + \sum \limits_{j=1}^{n}
\epsilon_{n,j}  b_{n-j} a_{j+1} + \epsilon_{n+2,0} a_{n+2} = 0,
\quad n = 1, 2, ... \ ,
$$
where $\epsilon_{n,j}$ are the entire powers of $e^{i \theta}$. This
successively implies the representations of $a_n$ by $b_j$ via
 \be\label{21}
a_n = (- 1)^{n-1} \epsilon_{n-1,0}  b_0^{n-1} - (- 1)^{n-1} (n - 2)
\epsilon_{1,n-3} b_1 b_0^{n-3} + \text{lower terms with respect to} \ b_0.
\end{equation}
This transforms any coefficient functional $J(w)$ on $S$ depending on a finite set of  distinguished
coefficients $a_{m_1}, \dots, a_{m_s}$  into a coefficient functional $\wt J(W^\mu)$
on $\Sigma^0$ depending on the corresponding coefficients $b_j$ and lifts the functionals  $\wt J(W)$ and $J(f)$ holomorphically onto the universal Teichm\"{u}ller space $\T$ modelled via bounded domain in the space $\B_2$ of Schwarzians $S_W^\mu$. In our case, $J(w) = c_n$ (expressed in terms of the corresponding coefficients $a_j$ of $w \in S$.

To lift $J$ onto the covering space $\T_1$, we again pass to functional $\wh J(\mu) = \wt J(W^\mu)$
lifting $J$ onto the ball $\Belt(\D)_1$ and apply the $\T_1$-equivalence, i.e., the quotient map
$$
\phi_{\T_1}: \ \Belt(\D)_1 \to \T_1, \quad \mu \to [\mu]_{\T_1}.
$$
Thereby the functional $\wt J(F^\mu)$ is pushed down to a bounded holomorphic functional $\mathcal J$
on the space $\T_1$ with the same range domain.

Regarding (after applying the Bers isomorphism theorem), the points of $\T_1$ as the pairs
$X_{F^\mu} = (S_{F^\mu}, F^\mu(0))$, where $\mu \in \Belt(\D)_1$ obey $\T_1$-equivalence, we come to a holomorphic functional
$$
\mathcal J(X_{F^\mu}) = \mathcal J(S_{F^\mu}, \ t), \quad t = F^\mu(0).
$$
on $\T_1 = \Fib (\T)$ and have to investigate its restriction to the image in $\Fib (\T)$ of the original submanifold $G$.

We shall need the following covering lemma generalizing Koebe's one-quarter theorem.
Let $\chi$ be a holomorphic map of a submanifold $G$ from of a complex Banach space $X$ into the space $\T$.

\noindent
{\bf Lemma 2}. \cite{Kr5} {\it Let $w(z)$ be a holomorphic univalent solution of the Schwarz differential equation $S_w(z) = \chi(\x)$
on $\D$ satisfying $w(0) = 0, \ w^\prime(0) = e^{i \theta}$ with the fixed $\theta \in [-\pi, \pi]$
and $\x \in G$ (hence $w(z) = e^{i \theta} z  + \sum_2^\infty a_n z^n$). Assume that
$$
a_{2,\theta}^0 = \sup \{|a_{2,\theta}|: \ S_w \in \chi(G)\}\ne 0,
$$
and let $w_0(z) = e^{i \theta} z + a_2^0 z^2 + \dots$ be one of the maximizing functions for $a_{2,\theta}^0$. Then:

(i) For every indicated function $w(z)$ , the image domain $w(D)$ covers entirely the disk
$D_{1/(2 |a_{2,\theta}^0|)} = \{|w| < 1/(2 |a_{2,\theta}^0|)\}$.

The radius value $1/(2 |a_{2, \theta}^0|)$ is sharp for this collection of
functions and fixed $\theta$, and the circle $\{|w| = 1/(2 |a_{2,\theta}^0|)$ contains points
not belonging to $w(\D)$ if and only if $|a_2| = |a_{2,\theta}^0|$
(i.e., when $w$ is one of the maximizing functions).

(ii) The inverted functions
$$
W(\zeta) = 1/w(1/\zeta) = e^{i \theta}\zeta - a_2^0 + b_1 \zeta^{-1} + b_2 \zeta^{-2} + \dots
$$
with $\z \in D^{-1}$ map domain $D^{-1}$ onto a domain whose boundary is entirely
contained in the disk $\{|W + a_{2,\theta}^0| \le |a_{2,\theta}^0|\}$.
}

This lemma yields that the boundary of domains $W^\mu(\D^*)$ for any $W^\mu(z) \in \Sigma^0$ are located in the disk $\{|W - b_0| \le |a_2^0| \}$  with $|a_{2, \theta}^0|$, and the second coordinate $t$ runs over some subdomain $D_\theta$ in the disk $\D_4 = \{|t| < 4\}$ containing the origin.
Since the functional $J$ is rotationally invariant, this subdomain $D_\theta$ is a disk $\D_{\a_\theta}$ of some radius $\a_\theta \le 2 |a_{2,\theta}^0|$.

We define on this disk the function
 \be\label{22}
\wt u_1(t) = \sup_{S_{W^\mu}} \mathcal J(S_{W^\mu}, t).
\end{equation}
taking the supremum over all $S_{W^\mu} \in \T$ admissible for a given $t = W^\mu(0) \in D_{\a_\theta}$, that means over the pairs $(S_{W^\mu}, t) \in \Fib(\T)$ with a fixed $t$ and pass to the upper semicontinuous regularization
$$
u_1(t) = \limsup\limits_{t^\prime \to t} \wt u_1(t^\prime).
$$

\noindent
{\it (b)} \ Now the crucial step in the proof of Theorem 2 is to establish that the function (22) inherits subharmonicity. We select on the unit circle a dense subset
$$
\mathbf e = \{z_1, z_2, \dots, z_n, \dots\}, \quad z_1 = e^{i \theta_1},
$$
and repeat successively for the above construction with fix points $z_1, z_2, \dots$ obtaining  similar
to (22) the corresponding functions $u_1(t), u_2(t), \dots$.
Let $u(t)$ be their upper envelope $\sup_n u_n(t)$ followed by its upper semicontinuous reglarization.

This construction is a special case of the general result considered in \cite{Kr5}, and this result
implies

\bigskip\noindent
{\bf Lemma 3}. {\it The functions $u(t)$ is logarithmically subharmonic in some disk $\D_\a$ with
$\a \le 2|a_2^0|$.  }

\bigskip\noindent
{\it (c)} \ Finally, we have to establish the range domain of $F^\mu(0)$ for $S_{F^\mu}$ running over
$\chi(G)$ and describe the boundary points of this domain.

The assumptions on $\chi(G)$ and the features of construction of the Bers fiber space
$\Fib(\T)$ imply that image of $\chi(G)$ in the space $\Fib(\T)$ also is a connected submanifold
covering $\chi(G)$. This allows one to apply to this image the same arguments as in \cite{Kr5}, \cite{Kr6}, applied there to the whole space $\Fib(\T)$.

It follows from the previous step that the indicated domain is rotationally symmetric and connected,
hence a disk $\D_\a = \{|t| < \a\}$ of some radius $\a \le 2 |a_2^0|$. Theorem 1 (parts {\it (a)} and {\it (b)}) yields that the extremal value of $\a$ equals $2|a_2^0|$.
The maximum of $|J|$ must  be attained on the boundary circle $\{|t| = 2 |a_2^0|\}$, corresponding by
Theorem 1 the function $w_0$, and the assertion of  Theorem 2 follows.

\bigskip\bigskip
\centerline{\bf 7. PROOF OF THEOREM 4}

\bigskip
To apply the above arguments, we need the following underlying lemma on openness
of the set of nonvanishing functions.

\bigskip\noindent
{\bf Lemma 4}. {\it Each point $f \in \mathcal B(H)$ has a neighborhood $U(f, \epsilon)$ in
$\mathcal B(H)$ filled by the functions which are zero free in the disk $\D$.
Take the maximal neighborhoods $U(f, \epsilon)$ with such property.
Then their union
$$
\mathcal U = \bigcup_{f\in \mathcal B(H)} U(f, \epsilon)
$$
is an open path-wise connective set, hence a domain, in the space $H$.       }

\bigskip
Observe that the assumption (11) plays also here a crucial role. It implies that in fact we have to establish the openness and connectivity of the set $\mathcal U$ in the $L_\iy$-norm.

\bigskip\noindent
{\bf Proof}. {\it (a) Openness}. It suffices to show that for each
$r > 1$ and $1 \le r^\prime < r $, every $f$ in the corresponding ball $\mathcal B(H(\D_r))$
of holomorphic functions on the larger disk $\D_r = \{|z| < r\}$
has a neighborhood $U(f, \epsilon)$ in $\mathcal B(H)$, which contains only
the functions that are zero free on the disk $\D$.

This is trivial for $ r > r^\prime > 1$. Let $r^\prime = 1$, and assume
the contrary. Then for some $r > 1$ there exist a function $f_0 \in \mathcal B(H(\D_r))$, a sequence of
functions  $f_n \in \mathcal B(H(\D_{r_n}))$ convergent to $f_0$ as $r_n \nearrow r$ so that
 \be\label{23}
\lim\limits_{n\to \iy} \|f_n - f_0\|_{H(\D_{r_n})} = 0,
\end{equation}
and a sequence of points $z_n \in \D_{r_n}$ convergent to $z_0$ with $|z_0| \le 1$ such that
$f_n (z_n) = 0 \ (n = 1, 2, \dots)$.

In the case $|z_0| < 1$, we immediately reach a contradiction, because then the uniform convergence of $f_n$ on compact sets in $\D$ implies $f_0(z_0) = 0$, which is impossible.

The case $|z_0| = 1$ requires other arguments. Since $f_0$ is holomorphic and does not vanish in the disk $\D_r$ with $r > 1$,
$$
\min_{|z| \le 1} |f_0(z)| = a > 0.
$$
Hence, for each $z_n$,
$$
|f_n(z_n) - f_0(z_n)| = |f_0(z_n)| \ge a,
$$
and by continuity, there exists a neighborhood $\Delta(z_n, \dl_n) =
\{|z - z_n| < \dl_n\}$ of $z_n$ in $\D$, in which
$$
|f_n(z) - f_0(z)| \ge a/2 \quad \text{for all} \ \  z.
$$
This yields
$$
\|f_n - f_0\|_{\B_2(D_{r_n})} \ge \|f_n - f_0\|_{\B_2} \ge  C(a),
$$
with a positive constant $C(a)$ depending only on $a$. Then the assumption (11) implies
$$
\|f_n - f_0\|_H \ge C_1(a),
$$
also for all $n$, but this contradicts to (23).

\bigskip\noindent
{\it (b)}  The path-wise connectedness of the set $\mathcal U$ is a simple consequence of
the fact that every its point $f(z) = c_0^1 + c_1^1 z + \dots$ is connected via homotopy
$f(z, t) = c_0^{0,j} + c_1^{0,j} t z + \dots, 0 \le t \le 1$, with the origin of $H$.
In view of the well-known properties of quasiconformal maps, this homotopy extends to
a complex holomorphic homotpy $\D \times \D \to B_2$, and then by (11) to homotopy
$\D \times \D \to H$.  The lemma follows.

Now Theorem 4 follows from Theorem 2 as a special case when $G = \mathcal B(H)$.

\bigskip
Note that openness given by Lemma 4 is also fulfilled for arbitrary Banach spaces obeying the assumption
(11).

\bigskip\bigskip
\centerline{\bf 8. GENERAL HILBERT SPACES OVER DISK}

\bigskip
We conclude with the following remarks.

\bigskip\noindent
{\bf 8.1}.
Theorem 2 also implies some useful consequences which concern solving the Hummel-Scheinberg-Zalcman problem on the general Hilbert spaces over the unit disk.

Namely, consider the union $\B_\iy = \bigcup_{p \ge 2} \B_p$ equipped with topology of inductive limit.
It contains the elements of any Hilbert space $H$ formed by holomorphic functions on $\D$ and their
limits in topology generated by the locally uniform convergence on $\D$.

Combining Theorems 3 and 4 and using this convergence, one obtains for such spaces the approximate values of the desired quantity (18).

\bigskip\noindent
{\bf 8.2}. An interesting {\bf open problem} is to find the extent in which the variational Theorem 1
is extended to more general Banach spaces.

Any such extension of this theorem gives rise to the corresponding extensions of Theorems 2 and 4.

\bigskip
{\small\it{ \leftline{Department of Mathematics, Bar-Ilan University, 5290002 Ramat-Gan, Israel} \leftline{and Department of Mathematics, University of Virginia,  Charlottesville, VA 22904-4137, USA}}


\bigskip
\bigskip
\begin{thebibliography}{EKK}
{\small

\bibitem{AW}
L.V. Ahlfors and G. Weill, {\it A uniqueness theorem for Beltrami equations}, Proc. Amer. Math. Soc. \textbf{13} (1962), 975-978.

\bibitem{Be}
L. Bers, {\it Fiber spaces over Teichm\"{u}ller spaces}, Acta Math. \textbf{130} (1973), 89-126.

\bibitem{EKK}
C.J. Earle, I. Kra and S. L. Krushkal, {\it Holomorphic motions and Teichm\"{u}ller spaces},
Trans. Amer. Math. Soc.  \textbf{343} (1994), 927-948.

\bibitem{GL}
F.P. Gardiner and N. Lakic, {\it Quasiconformal \Te \ Theory}, Amer. Math. Soc.,
Providence, RI, 2000.

\bibitem{HSZ}
J. A. Hummel, S. Scheinberg and L. Zalcman, {\it A coefficient problem for bounded nonvanishing functions}, J. Anal. Math. \textbf{31} (1977), 169-190.

\bibitem{Kr1}
S.L. Krushkal, {\it Quasiconformal Mappings and Riemann Surfaces}, Wiley, New York, 1979.

\bibitem{Kr2}
S.L. Krushkal, {\it Quasiconformal maps decreasing $L_p$ norm}, Siberian Math. J. \textbf{41} (2000), 884-888.

\bibitem{Kr3}
S.L. Krushkal, {\it Teichm\"{u}ller spaces and coefficient problems for univalent holomorphic functions},  Analysis and Mathematical Physics \textbf{10} (2020), no. 4.
https://doi.org/10.1007/s13324-020-00395-y

\bibitem{Kr4}
S.L. Krushkal, {\it Extremal quasiconformality vs rational approximation},
J. Math. Sci. \textbf{244} (1) (2020), 22-35; Ukr. Math. Bull. \textbf{16} (2019), 181-199.
DOI: 10.1007/s10958-019-04601-6

\bibitem{Kr5}
 S.L. Krushkal, {\it Towards to general distortion theory for univalent functions:
Teichm\"{u}ller spaces and coefficient problems of complex analysis}, J. Math. Sciences \textbf{299} (1)  (2026), 64-95; Ukr. Math. Bull. \textbf{22} (4) (2025), 548-590.
DOI 10.1007/s10958-026-08398-z

\bibitem{Kr6}
S.L. Krushkal, {\it A new method for estimating the coefficients of holomorphic functions},
Axioms \textbf{15}(5) (2026), 361; \ https://doi.org/10.3390/axioms15050361

\bibitem{Kz}
J. Krzyz, {\it Coefficient problem for bounded nonvanishing functions}, Ann. Polon. Math. \textbf{70} (1968), 314.

\bibitem{Le}
O. Lehto, {\it Univalent Functions and Teichm\"uller Spaces}, Springer, New York, 1987.

\bibitem{Ro}
H.L. Royden,
{\it Automorphisms and isometries of Teichmüller space}, Advances in the Theory of Riemann Surfaces (Ann. of Math. Stud., vol. 66), Princeton Univ. Press, Princeton, 1971, pp. 369-383.


}

\end{thebibliography}
\end{document}